\documentclass{elsart}
\usepackage{amssymb}
\usepackage{amsmath}
\journal{Journal of Pure and Applied Algebra}
\date{}

\newcommand\Qmax{Q^r_{\mathrm{max}}}
\newcommand\Qcl{Q^r_{\mathrm{cl}}}
\newcommand\f{{\mathcal F}}
\newcommand\te{{\mathcal T}}
\newcommand\ef{{\mathfrak F}}
\newcommand\vng{{\mathcal N}G}
\newcommand\ann{{\mathrm{ann}}}

\begin{document}
\begin{frontmatter}

\title{Differentiability of Torsion Theories}

\author{Lia Va\v s}

\address{Department of Mathematics, Physics and Computer Science,
University of the Sciences in Philadelphia, 600 S. 43rd St.,
Philadelphia, PA 19104}

\ead{l.vas@usip.edu}

\begin{abstract}
We prove that every perfect torsion theory for a ring $R$ is differential (in the sense of \cite{Bland_paper}). In this case, we construct the extension of a derivation of a right $R$-module $M$ to a derivation of the module of quotients of $M$. Then, we prove that the Lambek and Goldie torsion theories for any $R$ are differential.
\end{abstract}


\end{frontmatter}

\section{Introduction}

The many important examples of rings in analysis and differential geometry have inspired an interest to study rings equipped with maps that have properties of a derivation. As the derivation on a ring is not intrinsically a ring theoretic notion, it is of interest to study how a derivation agrees with other ring theoretic notions. In \cite{Golan_paper} and \cite{Bland_paper}, the authors study how derivations agree with an arbitrary hereditary torsion theory for that ring. Here, we continue this study by concentrating on some important classes and examples of hereditary torsion theories and prove that they are differential.

In Section \ref{basic_notions}, we recall
the definition of derivation on a ring, a (hereditary) torsion theory, Gabriel filter and rings and modules of quotients. We also recall some important examples of torsion theories (Lambek, Goldie, classical, etc.).

In Section \ref{section_perfect}, we study perfect torsion theories. These are torsion theories in which the module of quotients of every module is (isomorphic to) the tensor product of the module with the right ring of quotients. Thus, a perfect torsion theory is a generalization of the classical torsion theory of a right Ore ring. As the classical torsion theory is differential, this leads us to suspect that a perfect torsion theory might be differential too. The main result of Section \ref{section_perfect} states that this is indeed the case. Also, we give an explicit construction of the extension of derivation of a module to the derivation of the module of quotients of that module. We derive some corollaries of the main result.

In Section \ref{section_Lambek_and_Goldie}, we prove that the Lambek and Goldie torsion theories for every ring are differential.

\section{Torsion Theories, Differential Filters, Modules of Quotients}
\label{basic_notions}

Throughout this paper, $R$ denotes an associative ring with a unit. By a
module we mean a right module unless otherwise specified. We recall some definitions first.

A {\em derivation} on $R$ is a mapping $\delta: R \rightarrow R$ such that $\delta(r+s)=\delta(r)+\delta(s)$ and $\delta(rs)=\delta(r)s+r\delta(s)$ for all $r,s\in R.$ A mapping $d: M\rightarrow M$ on a right $R$-module $M$ is a {\em $\delta$-derivation} if
$d(x+y)=d(x)+d(y)$ and $d(xr)=d(x)r+x\delta(r)$ for all $x\in M$ and $r\in R.$ If $M$ is a left $R$-module, the notion of $\delta$-derivation on $M$ is defined analogously. The well known example of a derivation on a ring is obtained when considering $\frac{d}{dx}$ on the polynomial ring $R[x]$ over a ring $R.$ An {\em inner derivation} is another example of a ring derivation. If $R$ is a ring and $a\in R,$ define the derivation $\delta_a$ by $\delta_a(b)=ab-ba.$ If $M$ is an $R$-bimodule, then the map $d_a$ on $M$ given by $d_a(m)=am-ma$ is a $\delta_a$-derivation on $M.$ It is easy to see that the derivation $\frac{d}{dx}$ on $R[x]$ is not an inner derivation for some rings $R$ (e.g. take $R$ to be $\Qset$).

A {\em torsion theory} for $R$ is a pair $\tau = (\te, \f)$ of
classes of $R$-modules such that $\te$ and $\f$ are maximal classes having the property that
Hom$_R(T,F)=0,$ for all $T \in \te$ and $F \in \f.$
The modules in $\te$ are called {\em torsion modules} for $\tau$ and the modules in $\f$ are called
{\em torsion-free modules} for $\tau$.

A given class $\te$ is a torsion class of a torsion theory if an only if it is closed
under quotients, direct sums and extensions. A class $\f$ is a torsion-free class of a torsion theory if it is
closed under taking submodules, isomorphic images, direct products and extensions (see Proposition 1.1.9 in \cite{Bland_book}).

A torsion theory $\tau = (\te, \f)$ is {\em hereditary} if the
class $\te$ is closed under taking submodules (equivalently torsion-free class is closed under
formation of injective envelopes, see Proposition 1.1.6,
\cite{Bland_book}). The largest torsion theory in which a given class of
injective modules is torsion-free (the torsion theory {\em cogenerated} by that class) is hereditary. Some authors
(e.g. \cite{Golan_book}, \cite{Lambek}) consider just hereditary torsion theories.

A torsion theory $\tau = (\te, \f)$ is {\em faithful}
if $R\in \f.$

For every module $M$, the largest submodule of $M$ that belongs to $\te$ is called the
{\em torsion submodule} of $M$ and is denoted by $\te M$ (see Proposition 1.1.4 in \cite{Bland_book}). The
quotient $M/\te M$ is called the {\em torsion-free quotient} and
is denoted by $\f M.$

If $M$ is a right $R$-module with submodule $N$ and $m\in M,$ denote $\{r\in R\; | \; mr\in N\}$ by $(N : m).$ Then $(0:m)$ is the annihilator $\ann(m).$
A {\em Gabriel filter (or Gabriel topology)} $\ef$ on a ring $R$
is a nonempty collection of right ideals such that
\begin{enumerate}
\item If $I\in \ef$ and $r\in R,$ then $(I:r)\in \ef.$

\item If $I\in \ef$ and $J$ is a right ideal with $(J:r)\in
\ef$ for all $r\in I,$ then $J\in \ef$.
\end{enumerate}

If $\tau=(\te, \f)$ is a hereditary torsion theory, the collection of right
ideals
$\{ I\; |\; R/I\in \te\;\}$ is a Gabriel
filter. Conversely, if $\ef$ is a Gabriel
filter, then the class of modules $\{ M\; |\; \ann(m)\in \ef$ for every $m\in M\}$ is a torsion class of a
hereditary torsion theory. The details can
be found in \cite{Bland_book} or \cite{Stenstrom}.

We recall some important examples of torsion theories.
\begin{exmp}
{\em

(1) The torsion theory cogenerated by the injective
envelope $E(R)$ of $R$ is called the {\em Lambek torsion theory}.
It is hereditary, as it is cogenerated by an injective module, and
faithful. Moreover, it is the largest hereditary faithful torsion
theory. The Gabriel filter of this torsion theory is the set of
all dense right ideals (see Proposition VI 5.5, p. 147 in
\cite{Stenstrom}).

(2) The class of nonsingular modules over a ring $R$ is closed
under submodules, extensions, products and injective envelopes.
Thus, it is a torsion-free class of a hereditary torsion theory.
This torsion theory is called the {\em Goldie torsion theory}. It
is larger than any hereditary faithful torsion theory (see Example
3, p. 26 in \cite{Bland_book}). So, the Lambek torsion theory is
smaller than the Goldie's. If $R$ is right nonsingular, the Lambek and Goldie torsion theories
coincide (see  p. 26 in \cite{Bland_book} or p. 149 in \cite{Stenstrom}).

(3) If $R$ is a right Ore ring with the set of regular elements $T$ (i.e.,
$rT \cap tR \neq 0,$ for every $t \in T$ and $r\in R$), we can
define a hereditary torsion theory by the condition that a right
$R$-module $M$ is a torsion module iff for every $m\in M$, there
is a nonzero $t\in T$ such that $mt =0.$ This torsion theory is
called the {\em classical torsion theory of a right Ore ring}. It
is hereditary and faithful.

(4) Let $R$ be a subring of a ring $S$. The collection of all $R$-modules $M$ such that $M\otimes_R S = 0$
is closed under quotients, extensions and direct
sums. Moreover, if $S$ is flat as a left $R$-module, then this
collection is closed under submodules and, hence, defines a
hereditary torsion theory. We denote this torsion
theory by $\tau_S.$ As
all flat modules are $\tau_S$-torsion-free, $\tau_S$ is faithful.
If $R$ is a right Ore ring, then $\tau_{\Qcl(R)}$ is the classical
torsion theory.
}
\label{Examples}
\end{exmp}

A Gabriel filter $\ef$ is a {\em differential filter} if for every $I\in \ef$ there is $J\in \ef$ such that $\delta(J)\subseteq I$ for all derivations $\delta.$ The hereditary torsion theory determined by $\ef$ is said to be a {\em differential torsion theory} in this case. By Lemma 1.5 from \cite{Bland_paper}, a torsion theory is differential if and only if
\[d(\te M)\subseteq \te M\]
for every right $R$-module $M,$ every derivation $\delta$ and every $\delta$-derivation $d$ on $M.$

Examples 1.1 -- 1.4, page 3, from \cite{Bland_paper} provide some examples of differential torsion theories:
\begin{enumerate}
 \item  If $R$ is commutative, every Gabriel filter $\ef$ is differential as $\delta(I^2)\subseteq I$ for all $I\in \ef.$
\item If the torsion class $\te$ is closed for products, then it is differential (for details see Example 1.2, page 3 in \cite{Bland_paper}). As a consequence, every hereditary torsion theory for a left perfect ring is differential.
\item By Example 1.4, page 3, from \cite{Bland_paper}, the classical torsion theory of a right Ore ring is differential.
\end{enumerate}

If $\tau$ is a hereditary torsion theory with Gabriel filter $\ef$ and $M$ is a right $R$-module, define the module $M_{\ef}$
as the largest submodule $N$ of $E(M/\te M)$ such that $N/(M/\te M)$ is torsion module (i.e. the closure
of $M/\te M$ in $E(M/\te M)$). In \cite{Lambek} (definition of $Q(M)$ on page 11 and Example 5 on page 25) and \cite{Stenstrom} (exposition on pages 195--197) it is shown that $R_{\ef}$ has a
ring structure and that $M_{\ef}$ has a structure of a right
$R_{\ef}$-module.  The ring $R_{\ef}$ is called the {\em
right ring of quotients with respect to the torsion theory
$\tau$} and $M_{\ef}$ is called the {\em module
of quotients} of $M$ with respect to $\tau.$

For example, $\Qmax(R)$ is the right ring
of quotients with respect to the Lambek torsion theory (Example 1, page 200, \cite{Stenstrom}). $E(\f R)$ is the right ring of quotients with respect to the Goldie torsion theory if $\f R$ is the torsion free quotient of $R$ in the Goldie torsion theory (Propositions IX 1.7, 2.5, 2.7, and 2.11 and Lemma IX 2.10 in \cite{Stenstrom}). $\Qcl(R)$ is the right ring of quotients with respect to classical torsion theory if $R$ is right Ore (Example 2, page 200, \cite{Stenstrom}).

Consider the map $\phi_M:M\rightarrow M_{\ef}$ obtained by
composing the projection $M\rightarrow M/\te M$ with the injection $M/\te M\rightarrow M_{\ef}.$ Then,
the kernel and cokernel of $\phi_M$ are torsion modules and $M_{\ef}$ is torsion-free (Lemmas 1.2 and 1.5, page 196, in \cite{Stenstrom}). In Corollary 1 of \cite{Golan_paper}, Golan has shown that for the differential filter $\ef$ and the $\delta$-derivation $d_M$ of a right $R$-module $M$, $d_M$ extends to a derivation $d_{M_{\ef}}$ of the module of quotients $M_{\ef}$ such that $d_{M_{\ef}}\phi_M=\phi_M d_M.$ Bland proved that such extension is unique and that the converse is also true. Namely, Proposition 2.3 of his paper \cite{Bland_paper} states the following.

\begin{thm} Let $\ef$ be a Gabriel filter. The filter $\ef$ is differential if and only if
every derivation on any module $M$ extends uniquely to a derivation on the module of quotients $M_{\ef}.$
\label{Blands_Theorem}
\end{thm}

\section{Differentiability of a Perfect Torsion Theory}
\label{section_perfect}

A ring homomorphism $f:R\rightarrow S$ is
a {\em ring epimorphism} if for all rings $T$ and
homomorphisms $g,h: S\rightarrow T,$ $gf = hf$ implies $g=h.$
The situation when $S$ is flat as left $R$-module is of special
interest. If $f:R\rightarrow S$ is a ring epimorphism with $S$ flat as left $R$-module,
then $S$ is a right ring of quotients with respect to the torsion theory with Gabriel filter
$\ef=\{I | f(I)S=S\}.$
$S$ is called a {\em perfect right ring of quotients,} a {\em flat epimorphic extension} of $R,$ or
a {\em perfect right localization of $R$} in this case. A hereditary torsion theory $\tau$ with Gabriel filter $\ef$ is
called {\em perfect} if the right ring of quotients $R_{\ef}$ is
perfect and $\ef=\{I| \phi_I(I)R_{\ef}=R_{\ef}\}$. The Gabriel filter
$\ef$ is called {\em perfect} in this case.

Perfect filters have a nice description. For a Gabriel filter
$\ef,$ let us look at the canonical maps $i_M: M\rightarrow
M\otimes_R R_{\ef}$ and $\phi_M: M\rightarrow M_{\ef}.$ There is a
unique $R_{\ef}$-map $\Phi_M: M\otimes_R R_{\ef}\rightarrow
M_{\ef}$ such that $\phi_M = \Phi_M i_M.$ The perfect filters are
characterized by the property that the map $\Phi_M$ is an isomorphism
for every module $M.$ More details can be found in \cite{Stenstrom} (Theorem XI 3.4, p.
231). Also, a perfect filter coincides with the filter of torsion theory $\tau_{R_{\ef}}$ obtained by tensoring
with the ring of quotients $R_{\ef}$ (for details see Lemma 8 in \cite{Lia4}).

We show that all perfect filters are differential. First we need an easy lemma.
\begin{lem}
If $M$ is a right $R$-module, $N$ is an $R$-bimodule, $d_M$ and $d_N$ $\delta$-derivations and $1_M$ and $1_N$ identity mappings on $M$ and $N$, then the map $d=d_{M\otimes_R N}: M\otimes_R N\rightarrow  M\otimes_R N$ defined by $d=d_M\otimes 1_N +1_M\otimes d_N$ is a $\delta$-derivation.
\label{derivation_of_tensor}
\end{lem}
\begin{pf} Clearly the map $d$ is additive.
Let $m\in M,n\in N$ and $r\in R$ be arbitrary. $d((m\otimes n)r)=d(m\otimes nr)=d_M(m)\otimes nr + m\otimes d_N(nr)=
d_M(m)\otimes nr + m\otimes d_N(n)r+ m\otimes n\delta(r)=(d(m\otimes n))r+(m\otimes n)\delta(r).$ From this observation and the additivity of $d$ it follows that $d$ is a $\delta$-derivation. \qed
\end{pf}

Now we show the main result of this section.
\begin{prop}
If a Gabriel filter $\ef$ is perfect, then it is differential.
\label{perfect_is_differential}
\end{prop}
\begin{pf}
By Theorem \ref{Blands_Theorem}, it is sufficient to show that the derivation $d_M: M\rightarrow M$ extends uniquely to a derivation $d_{M_{\ef}}: M_{\ef}\rightarrow M_{\ef}$ such that $d_{M_{\ef}}\phi_M=\phi_M d_M,$ for every module $M.$ This is automatically fulfilled if $M$ is torsion-free by Corollary 2.2 in \cite{Bland_paper}. If $\ef$ is perfect, $R$ is torsion-free (as the torsion submodule of $R$ with respect to a perfect torsion theory is isomorphic to Tor$^R_1(R,R_{\ef}/R)=0$, see Theorem XI 3.4, p.
231 in \cite{Stenstrom}) so we obtain the unique extension $\delta_{R_{\ef}}$ of $\delta$ to a derivation of $R_{\ef}.$

As $\ef$ is perfect, the unique map $\Phi_M: M\otimes_R R_{\ef}\rightarrow
M_{\ef}$ such that $\phi_M = \Phi_M i_M$ is an isomorphism for every module $M.$
Define \[d_{M_{\ef}}=\Phi_M d_{M\otimes_R R_{\ef}}\Phi_M^{-1}\]
where map $d_{M\otimes_R R_{\ef}}$ is map from Lemma \ref{derivation_of_tensor} defined via $d_M$ and $\delta_{R_{\ef}}.$

Clearly, $d_{M_{\ef}}$ is additive and a short, straightforward calculation shows that it is a derivation. We show that the following diagram commutes
\[ \begin{array}{lclcl}
M & \stackrel{i_M}{\longrightarrow} & M\otimes_R R_{\ef} & \stackrel{\Phi_M}{\longrightarrow} & M_{\ef}\\

\downarrow d_M & & \downarrow d_{M\otimes_R R_{\ef}} & & \downarrow d_{M_{\ef}} \\
M & \stackrel{i_M}{\longrightarrow} & M\otimes_R R_{\ef} & \stackrel{\Phi_M}{\longrightarrow} & M_{\ef}\\
\end{array} \]

As $\delta(1)=0,$ $\delta_{R_{\ef}}(1)=0.$  Thus, if $m\in M$ is arbitrary, $d_{M\otimes_R R_{\ef}}i_M(m)=d_{M\otimes_R R_{\ef}}(m\otimes 1)=d_{M}(m)\otimes 1+0=i_M(d_M(m)).$ So, $d_{M\otimes_R R_{\ef}}i_M=i_M d_M.$
$d_{M_{\ef}}\Phi_M=\Phi_M d_{M\otimes_R R_{\ef}}$ by definition of $d_{M_{\ef}}.$
This gives us

\[\begin{array}{rcll}
d_{M_{\ef}}\phi_M & = & \Phi_M d_{M\otimes_R R_{\ef}}\Phi_M^{-1}\phi_M & (\mbox{definition of }d_{M_{\ef}})\\
 & = & \Phi_M d_{M\otimes_R R_{\ef}}i_M  & (\mbox{as }\phi_M = \Phi_M i_M)\\
& = & \Phi_M i_M d_M & (\mbox{by the above diagram})\\
& = & \phi_M d_M  & (\mbox{as }\phi_M = \Phi_M i_M).
\end{array}\]
Finally, $d_{M_{\ef}}$ is unique by Proposition 2.1 in \cite{Bland_paper}. \qed
\end{pf}

\begin{rem} {\em
It is possible to give the explicit description of the derivation $\delta_{R_{\ef}}$ extending $\delta$ on the right ring of quotients $R_{\ef}$ if $\ef$ is perfect. Namely, if $\ef$ is perfect, then every element $q\in R_{\ef}$ has the
property
\[qr_j\in R\mbox{ and }\sum_{j=1}^n r_j q_j=1\mbox{ for some }n, q_j\in R_{\ef}\mbox{ and }r_j\in R,\;j=1,\ldots,
n.\] For proof see Theorem 2.1, p. 227 of \cite{Stenstrom}.
By the above property, $q=q 1 = \sum  q r_j q_j.$ To define
$\delta_{R_{\ef}}$ to be a derivation extending $\delta$ we need $\delta_{R_{\ef}}(q)=\sum (\delta(q r_j)q_j+ q r_j \delta_{R_{\ef}}(q_j)).$ To define this, it is sufficient to define the sum $\sum r_j \delta_{R_{\ef}}(q_j).$
As $\sum r_j q_j=1,$ we can define $\sum r_j \delta_{R_{\ef}}(q_j)$ to be $\delta(1)-\sum \delta(r_j)q_j=-\sum \delta(r_j)q_j.$ Thus define
\[\delta_{R_{\ef}}(q)=\sum_{j=1}^n (\delta(q r_j)q_j-q\delta(r_j)q_j)\]
}
\end{rem}

\begin{cor} If $\ef$ is a perfect filter, $M$ a right $R$-module and $d_M$ derivation on $M$, then the unique derivation extending $d_M$ to the module of quotients $M_{\ef}$ is given by
\[d_{M_{\ef}}(x)=\Phi_M\left(\sum_i d_M(m_i)\otimes q_i+m_i\otimes \sum_j (\delta(q_i r_{ij})q_{ij}-q_i\delta(r_{ij})q_{ij})\right)\]
where  $\Phi_M(\sum_i m_i\otimes q_i)=x,$ $m_i\in M$ and $q_i\in R_{\ef}$ and
$r_{ij}\in R, q_{ij}\in R_{\ef}$ are such that $q_i r_{ij}\in R$ for all $i$ and $j$ and
$\sum_j r_{ij}q_{ij}=1$ for all $i.$
\end{cor}
\begin{pf}
The proof follows  directly from the remark above and the proof of Proposition \ref{perfect_is_differential}. \qed
\end{pf}

The converse of the Proposition \ref{perfect_is_differential} is not true. Any infinite group $G$ gives rise to the group von Neumann algebra $\vng$ that is not
semisimple (see Exercise 9.11, p. 367 in \cite{Lu_book}) and so the Lambek torsion theory of $\vng$ is not perfect (by Theorem 23 (3) in \cite{Lia2}, Theorem 12 from \cite{Lia4} and Theorem 4 in \cite{Lia1}).
However, if $G$ is abelian and infinite, every torsion theory for $\vng$ is differential as $\vng$ is abelian. Another example of a ring with the Lambek torsion theory that is differential but not perfect is the ring from Example 4.2 in \cite{Lia4}.

The following is a direct corollary of Proposition \ref{perfect_is_differential}.
\begin{cor}
If $R$ is right noetherian and hereditary, then all hereditary torsion theories are differential. In particular, a semisimple ring has all hereditary torsion theories differential.
\end{cor}
\begin{pf}
If $R$ is right noetherian and hereditary, then all torsion theories are perfect (Corollary 3.6, p. 232 in \cite{Stenstrom}). \qed
\end{pf}

\begin{cor}
If $R$ is right semihereditary and $S$ is an extension of $R$ that is flat as left $R$-module, then $\tau_S$ is differential.
\end{cor}
\begin{pf}
We show that the assumptions of the claim imply that $\tau_S$ is perfect. Then the claim follows from Proposition \ref{perfect_is_differential}.

By the definition of $\tau_S,$ its filter $\ef$ consists of the right ideals $I$ such that $IS=S.$ This filter has a base of finitely generated ideals as if $I\in \ef, $ then $1\in S=IS$ so there is a nonnegative integer $n,$ $r_i\in I$ and $s_i\in S,$ $i=1,\ldots,n$ such that $\sum r_i s_i=1.$ Thus the ideal $J=\langle r_1,\ldots,r_n\rangle$ is finitely generated ideal in $\ef$ with $J\subseteq I.$

As $R$ is right semihereditary, every finitely generated right ideal is projective. Thus, $\ef$ has a base of projective ideals. But then $\tau_S$ is perfect by Proposition 3.3, p. 230 and Proposition 3.4, p. 231 in \cite{Stenstrom}. \qed
\end{pf}

\section{Differentiability of the Lambek and Goldie Torsion Theories}
\label{section_Lambek_and_Goldie}

First, we show that the Lambek torsion theory is differential. Recall that for a right module $M$ and $x\in M,$ $x$ is in the torsion submodule with respect to Lambek torsion theory if and only if $\ann(x)$ is a dense right ideal. Recall also that a right ideal $I$ is dense if and only if for every $r,s\in R$ such that $s\neq 0,$ there is $t\in R$ such that $st\neq 0$ and $rt\in I$.
\begin{prop}
The Lambek torsion theory for any ring $R$ is differential.
\label{Lambek_is_diff}
\end{prop}
\begin{pf}
Let $M$ be any $R$-module and $x$ any element of $M$. For the Lambek torsion theory to be differential it is sufficient to show that $\ann(d(x))$ is dense, whenever $\ann(x)$ is a dense right ideal.

Let $\ann(x)$ be dense and let $r,s\in R,$ $s\neq 0.$  As $\ann(x)$ is dense, there is $t_1\in R$ such that $st_1\neq 0$ and $xrt_1=0.$ As $\ann(x)$ is dense, for $\delta(rt_1)$ and $st_1\neq 0$ there is $t_2\in R$ such that $st_1 t_2\neq 0$ and $x\delta(rt_1)t_2=0.$

$xrt_1=0$ implies that $d(xrt_1)=0$. Thus $d(x)rt_1=-x\delta(rt_1).$ So, $d(x)rt_1 t_2=-x\delta(rt_1)t_2=0.$ Thus, for $r$ and $s\neq 0,$ we have found $t=t_1t_2$ such that $rt\in \ann(d(x))$ and $st\neq0.$ Hence, $\ann(d(x))$ is dense. \qed
\end{pf}

If $R$ is right nonsingular, the Goldie and Lambek torsion theories coincide. Thus, from Proposition \ref{Lambek_is_diff} it follows that the Goldie torsion theory of a right nonsingular ring is differential. We show that this is the case even if $R$ is not right nonsingular. Let $\ef_G$ denotes the Gabriel filter of the Goldie torsion theory. Then all the essential ideals are in $\ef_G$ and
\[\ef_G=\{I\; | \mbox{ there exists } J, I\subseteq J, J\subseteq_e R\mbox{ and }(I : j)\subseteq_e R\mbox{ for all }j\in J \}.\]
For proof see Proposition 6.3, p. 148 in \cite{Stenstrom}. From this observation it is easy to see that
\[\ef_G=\{I\; | \{r\in R\; |\; (I:r) \subseteq_e R\} \subseteq_e R \}.\]

Thus, if $M$ is a right $R$-module and $x$ is an element of $M,$ $x$ is in the torsion submodule for Goldie torsion theory (i.e. $\ann(x)\in \ef_G$) if and only if
\[\{r\in R\; |\; \ann(xr) \subseteq_e R\} \subseteq_e R.\] We shall follow the notation from \cite{Lam} (see p. 255 of \cite{Lam}) and denote \[
\begin{array}{rcl}
\ann(x)^* & = & \{r\in R\; |\; \ann(xr) \subseteq_e R\}\mbox{ and }\\
\ann(x)^{**} & = & \{r\in R\; |\; \{s\;|\; \ann(xrs) \subseteq_e R\}\subseteq_e R\}.
\end{array}
\]
Using this notation, $\ann(x)\in\ef_G$ if and only if $\ann(x)^*\subseteq_e R.$

To prove the main result, we need four lemmas.
\begin{lem}  If $M$ is a right $R$-module and $x\in M,$ then the following are equivalent
\begin{enumerate}
\item $\ann(x)\in \ef_G,$

\item $\ann(x)^*\subseteq_e R,$

\item $\ann(x)^{**}\subseteq_e R.$
\end{enumerate}
\label{claim_1}
\end{lem}
\begin{pf}
(1) and (2) are equivalent by the observation above. (2) implies (3)
as $\ann(x)^*\subseteq \ann(x)^{**}$ (for details see p.255 in \cite{Lam}). (3) implies (2) by the transitivity of relation $\subseteq_e$ and because $\ann(x)^*\subseteq_e \ann(x)^{**}$ (Proposition 7.29 (2), p. 255 in \cite{Lam}).
\qed
\end{pf}

\begin{lem} Let $M$ be a right $R$-module and $x,y\in M.$ Then
\begin{enumerate}
\item $\ann(x)\subseteq_e R$ implies $\ann(xr)\subseteq_e R$ for any $r\in R.$
\item $\ann(x)^*\subseteq_e R$ implies $\ann(xr)^*\subseteq_e R$ for any $r\in R.$
\item $\ann(x)^*\subseteq_e R$ and $\ann(y)^*\subseteq_e R$ imply $\ann(x+y)^*\subseteq_e R.$
\end{enumerate}
\label{claim_1.5}
\end{lem}
\begin{pf}
(1) follows from Lemma 7.2, p. 246 in \cite{Lam}. (2) and (3) are true as $x,y\in \te M$ implies that $xr\in \te M$ and $x+y\in \te M$ where $\te M$ is the torsion submodule of $M$ for the Goldie torsion theory.
\qed
\end{pf}

\begin{lem} If $M$ is a right $R$-module and $x\in M$ is such that $\ann(x)\subseteq_e R,$ then $\ann(d(x))^*\subseteq_e R.$
\label{claim_2}
\end{lem}
\begin{pf}
Let $r\in R.$ As $\ann(x)\subseteq_e R,$ there is $s\in R,$ $rs\neq 0$ and $xrs=0.$ Thus $0=d(xrs)=d(x)rs+x\delta(rs).$ As $\ann(x)$ is essential, $\ann(x\delta(rs))$ is also essential by part (1) of Lemma \ref{claim_1.5}. Thus, $\ann(d(x)rs)=\ann(-x\delta(rs))$ is essential as well and $0\neq rs\in\ann(d(x))^*.$ So
$\ann(d(x))^*$ is essential. \qed
\end{pf}

\begin{lem} If $M$ is a right $R$-module and $x\in M$ is such that $\ann(x)^*\subseteq_e R,$ then $\ann(d(x))^{**}\subseteq_e R.$
\label{claim_3}
\end{lem}
\begin{pf} Let $r$ be arbitrary element of $R$. As $\ann(x)^*\subseteq_e R,$ there is $s\in R,$ $rs\neq 0$ and $\ann(xrs)\subseteq_e R.$ By Lemma \ref{claim_2},
$\ann(d(xrs))^*\subseteq_e R.$ By part (2) of Lemma \ref{claim_1.5}, $\ann(x)^*\subseteq_e R$ implies that $\ann(x\delta(rs))^*\subseteq_e R.$ But then $\ann(d(x)rs)^*=\ann(d(xrs)-x\delta(rs))^*$ is essential in $R$ by part (3) of Lemma \ref{claim_1.5}.
Thus, $\{\;t\; |\; \ann(d(x)rst)\subseteq_e R\}\subseteq_e R$ and so $0\neq rs \in \ann(d(x))^{**}.$
\qed
\end{pf}

\begin{prop}
The Goldie torsion theory for any ring $R$ is differential.
\label{Goldie_is_diff}
\end{prop}
\begin{pf}
Let $M$ be a right $R$-module and $x\in M$ such that $\ann(x)\in \ef_G.$ Then $\ann(x)^*\subseteq_e R$ by Lemma \ref{claim_1} and so $\ann(d(x))^{**}\subseteq_e R$ by Lemma \ref{claim_3}. But this implies that $\ann(d(x))\in \ef_G$ by Lemma \ref{claim_1}. \qed
\end{pf}

Theorem \ref{Blands_Theorem}, Proposition \ref{perfect_is_differential}, Proposition \ref{Lambek_is_diff}, and
Proposition \ref{Goldie_is_diff} yield the following corollary.

\begin{cor}
If a Gabriel filter $\ef$ is perfect, Lambek or Goldie, then every derivation on any module $M$ lifts uniquely to a derivation on the module of quotients $M_{\ef}.$ In particular, every derivation on $R$ uniquely lifts to a derivation of $R_{\ef}.$
\end{cor}

In this paper, we have shown that some important examples of hereditary torsion theories are differential. This leads us to raise the question whether {\em every} hereditary torsion theory is differential.

\end{document}